\documentclass [12pt,a4paper,reqno]{amsart}
  \input amssymb.sty
\usepackage{amsmath}
\usepackage{amsfonts}
\usepackage{epsfig}

 \newtheorem{thm}{Theorem} [section]

\newtheorem*{thm*}{Theorem}

\newtheorem{lem} {Lemma}[section]
\newtheorem{lemma} {Lemma}[section]
\newtheorem{prop} {Proposition}[section]

\newtheorem{acknowledgment*}[thm] {Acknowledgment}

 \theoremstyle{definition}
 \newtheorem*{remark*}{Remark}
  \newtheorem{defn} {Definition}[section]
 \newtheorem{remark}{Remark}[section]
 \newtheorem*{defn*}{Definition}

\newtheorem*{notation*}{Notation}

\newcommand{\lemref}[1]{Lemma~\ref{#1}}

\newcommand{\Aut}{\operatorname{Aut}}
\newcommand{\codim}{\operatorname{codim}}
\newcommand{\Center}{\operatorname{Center}}
\newcommand{\Vertices}{\operatorname{Vertices}}

\newcommand{\p}{\pi_1}

\newcommand{\vp}{\varphi}
\newcommand{\Z}{\mathbb {Z}}
\newcommand{\C}{\mathbb {C}}

\renewcommand{\P}{\mathbb {P}}

  \newcommand{\B}{\mathcal{B}}
 \newcommand{\PP}{\mathcal{P}}

\begin{document}
\title [On fundamental groups related to  the
Hirzebruch surface $F_1$]{The fundamental group of the complement of
the branch curve of the Hirzebruch surface $F_1$}
\author[M. Friedman]{Michael Friedman}
\author[M. Teicher]{Mina Teicher}

\maketitle \numberwithin{equation}{section}
\maketitle
\begin{abstract}
Given a projective surface and a generic projection to the plane,
the braid monodromy factorization (and thus, the braid monodromy
type) of the complement of its branch curve is one of the most
important topological invariants (\cite{KuTe}), stable on
deformations. From this factorization, one can compute the
fundamental group of the complement of the branch curve, either in
$\C^2$ or in $\C\P^2$. In this article, we show that these groups,
for the Hirzebruch surface $F_{1,(a,b)}$, are almost-solvable. That
is - they are an extension of a solvable group, which strengthen the
conjecture
on degeneratable surfaces (see \cite{Te1}).\\\\
\textbf{keywords}:Hirzebruch surfaces, degeneration, generic
projection, branch
curve, braid monodromy, fundamental group, classification of surfaces.\\
\textbf{AMS classification numbers}: 14D05, 14D06, 14E25, 14H30,
14J10, 14Q05, 14Q10.
\end{abstract}

 \tableofcontents

\section{Introduction}

 In the study of smooth algebraic surfaces of degree $n$, which are embedded in
 $\mathbb C\mathbb P^N$, one can consider the surface $X$ as a branched cover
 of $\mathbb C\mathbb P^2.$ In this case the branch locus, $S_X$ in $\mathbb
 C\mathbb P^2,$ plays a crucial role. It is, in general, singular and, if the
 projection $X\to\C\P^2$ is generic, the singularities are nodes and cusps.
 The significance of $S_X$ (or of $S\subset \C^2\subset\C\P^2,$ a generic affine
 portion of $S_X$) arises when studying equivalence class of the
 braid monodromy factorization of the branch curve $S_X$ (which is known to be
 the BMT invariant of the surface $X$; see \cite{Te1}). From this
 factorization one can induce the fundamental groups $\overline
 G=\pi_1(\C\P^2-S_X)$ or $G=\pi_1(\C^2-S),$ which are stable on
 deformations. That is, if two surfaces have different fundamental groups, then
 they are not deformation equivalent.  For surfaces $X,\,Y$ denote
 $X \overset{G}{\backsimeq} Y\,\Leftrightarrow\, G_X = G_Y$ and
 $\overline G_X = \overline G_Y$; $X \overset{Diff}{\backsimeq} Y \Leftrightarrow$
 $X$ is diffeomorphic to $Y$; $X \overset{Def}{\backsimeq} Y \Leftrightarrow$
 $X$ is deformation equivalent to $Y$; and $X \overset{BMT}{\backsimeq} Y \Leftrightarrow$
 $X$ and $Y$ has the same BMT invariant.

 It turns out that $X \overset{Def}{\backsimeq}
 Y \Rightarrow X \overset{G}{\backsimeq} Y$ but the inverse
 direction is not correct (see \cite{KhKu}); and $X \overset{Def}{\backsimeq}
 Y \Rightarrow X \overset{BMT}{\backsimeq} Y \Rightarrow X \overset{Diff}{\backsimeq}
 Y$ (and again - the inverse directions are not correct; see \cite{CW},\cite{KhKu}).

 In this article, we take $X$ to be the Hirzebruch surface $F_1;$ this surface
 is the projectivization of the line bundle $\mathcal
 O_{\C\P^1}(1)\oplus\mathcal O_{\C\P^1}.$ We then embed it in $\C\P^n$ with
 respect to the linear system $|aC+bE_0|,$ where $C,E_0$ generate the Picard
 group of $F_1, \,b>1, a\geq1.$ We show that $G$ and $\overline G$ can be computed
 when $X=F_{1,(a,b)}$, which is the image of $F_1$ after the embedding w.r.t.
 the above linear system.

 It is conjectured (\cite{Te1}) that $G$ and $\overline G$ are almost solvable in a large family of
 surfaces: that is, these groups are
 extensions of a solvable group by the symmetric group. So far, it
 was proven for $V_p$ (the Veronese surface; \cite{Te}) and $X_{p,q}$ (the double-double
 covering of $\C\P^1\times \C\P^1$; \cite{Mo}).

Our main result proves that $X=F_{1,(a,b)}$\ $(b>1, a\geq1)$
satisfies the conjecture. In particular ,there exists a series
$$1\triangleleft A_1\triangleleft A_2\triangleleft A_3\triangleleft G$$
s.t.
$$G/A_3\simeq S_{2ab+b^2},$$
$$A_3/A_2\simeq \mathbb Z,$$
$$A_2/A_1\simeq(\mathbb Z_{b-2a})^{2ab+b^2-1}$$
$$A_1\simeq\begin{cases}\mathbb Z_2\quad & b\ \text{even},\ a\ \text{odd}\\
1\quad &\text{otherwise}\end{cases}$$ and a series
$$1\triangleleft \overline A_1\triangleleft \overline A_2\triangleleft \overline A_3\triangleleft
\overline G$$ where
$$\overline G/\overline A_3=G/A_3,$$
$$\overline A_3/\overline A_2\simeq \mathbb Z_m,\quad
m=3ab-a-b+\frac{3b^2-3b}{2},$$
$$\overline A_2/\overline A_1=A_2/A_1,$$
$$\overline A_1=A_1.$$

As noted, the significance of this article lies in the fact that $G$
and $\overline G$ are determined by the deformation type, since they
are stable under deformation of the surface. Thus, computing $G$ and
$\overline G$ explicitly (and the series of groups derived from
them) can help us distinguish between non--deformation equivalent
Hirzebruch surfaces.

Another important aspect of this article is the fact that it gives a
general approach and another example of how to compute and deal with
the fundamental groups $G$ and $\overline G$. So far, only a few
examples of calculating these groups were presented (see
\cite{BGT5}, \cite{Robb}), and most of the calculations dealt with
the Galois cover of such a degeneratable surface; especially with
finding the fundamental group of this Galois cover (see
\cite{MoRoTe}, \cite{MoTe}). Calculating $G$ and $\overline G$ is
another step in understanding the whole structure of these groups
with respect to surfaces which can be degenerated.

\section{Hirzebruch surfaces and their degenerations}

The Hirzebruch surfaces $F_k$ (for $k\ge0$) are given by the equation
$x_1y_1^k=x_2y_2^k$ in $\C\P^1\times\C\P^1$. However, the construction these
days is as follows: the $k$-th Hirzebruch surface is the projectivization of
the vector bundle $\mathcal O_{\C\P^1}(k)\oplus \mathcal O_{\C\P^1}.$

Let $\sigma$ be a holomorphic section of $\mathcal O_{\C\P^1}(k),$ and let
$E_0\subset F_k$ denote the image of the section $(\sigma,1)$ of $\mathcal
O_{\C\P^1}(k)\oplus\mathcal O_{\C\P^1}.$ The curve $E_0$ is called a
\textit{zero section} of $F_k.$ All zero sections are homologous and hence
define a divisor class which is independent of choice of $\sigma.$ Let $C$
denote a fiber of $F_k.$ The Picard group of $F_k$ is generated by $E_0$ and
$C.$ It is elementary that $E_0^2=k,$\ $C^2=0$ and $E_0\cdot C=1.$

The surface $F_0$ is the quadric $\C\P^1\times \C\P^1,$ and $F_1$ is
the blow-up of the plane $\C\P^2.$ For $k>0,$ the surface $F_k$
contains a unique (irreducible) curve of negative self-intersection
$-k.$ This curve is a section of the bundle; it is denoted
$E_\infty$ and is called the \textit{negative section} or the
\textit{section at infinity}. We mention that it can be contracted
to an isolated normal singularity, the resulting normal surface
being the cone over the rational normal curve of degree $k.$ Zero
sections are always disjoint to $E_\infty.$ Schematically, we
describe $F_k$ as in Fig.~1.1.\\

 \begin{center}
\epsfig{file=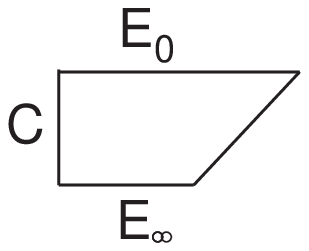}\\
(figure 1.1)
\end{center}

Let $F_k$ be the $k$-th Hirzebruch surface. Let $E_0,$ $E_\infty$, $C$ be as in
the Introduction. For $a,b\ge1$, or for $a=0 $ and $k\ge1,$ the divisor
$ac+bE_0$ on $F_k$ is very ample and thus defines an embedding $f_{|aC+bE_0|}:
F_k\hookrightarrow \C\P^N$. Let $F_{k(a,b)}=f_{|aC+bE_0|}(F_k)$\ $(\subseteq
\C\P^N).$ For  $k>0,$ the map $f_{|0\cdot C+bE_0|}$ collapses the section at
infinity to a point, so $F_{k(0,b)}$ is the image of the cone over the rational
normal curve of degree $k$ with respect to a suitable embedding.

In \cite{MoRoTe}, a degeneration to a union of $2ab+kb^2$ planes was
constructed in the following configuration (in Fig.~1.2, $k=2$,\
$a=2,$\ $b=3$ was taken). Each triangle represents a plane and each
inner edge represents an intersection line between planes.

 \begin{center}
\epsfig{file=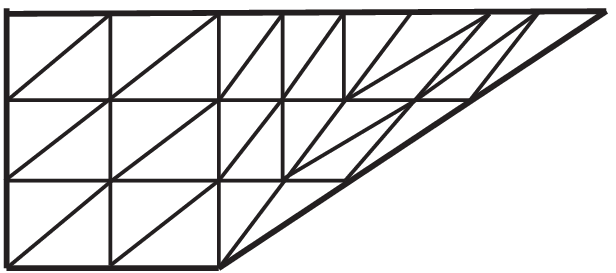}\\
(figure 1.2)
\end{center}
This degeneration is obtained using a technique developed by
Moishezon-Robb-Teicher which they refer to as the D-construction.
The D-construction is described (and prove to work) in \cite{BGT5}.
Specific degeneration for the Hirzebruch surfaces using the
D-construction is explained in \cite[Section 2, Theorem
2.1.2]{MoRoTe}. The difference between the D-construction and other
blow-up procedures for obtaining degenerations is that the
D-construction can also be applied along a subvariety of $\codim 0$
(see, for example, Step 2 below). The degeneration is obtained via
the following steps.
\begin{enumerate}
\item[1.] D-construction along $C$ to get $F_{0(a,b)}\cup F_{k(a-1,b)}.$
\item[2.] D-construction along $F_{0(1,b)}$ to get $F_{0(1,b)}\cup
F_{0(1,b)}\cup  F_{k(a-2,b)}.$ \item[3.] Induction on the second
step to get $\underbrace{F_{0(1,b)}\cup\dots\cup F_{0(1,b)}}_{a\
\text{times}}\cup F_{k(0,b)}$ (see \cite{MoTe5}). \item[4.]
Degeneration of each $F_{0(1,b)}$ to a union of $2b$ planes in the
following configuration (here  $b=3).$

 \begin{center}
\epsfig{file=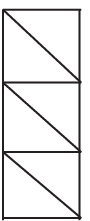}\\
(figure 1.3)
\end{center}
\item[5.] D-construction on $F_{k(0,b)}$ to get
$\underbrace{F_{1(0,b)}\cup\dots\cup F_{0(1,b)}\cup F_{1(0,b)}}_{k\
\text{times}}$
\newline $(F_{1(0,b)}$ is the Veronese
surface $V_b$). \item[6.] Degeneration of each $F_{1(0,b)}$ to a
union of $b^2$ planes in the following configuration (here $b=3):\\$

 \begin{center}
\epsfig{file=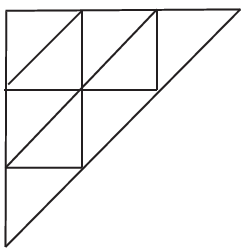}\\
(figure 1.4)
\end{center}
\end{enumerate}

Note that in our case $k=1$; so we are looking at the surface $F_{1(a,b)}.$ We
now describe in greater detail the degenerated object and its branch curve
using the degeneration described earlier. $F_{1(a,b)}$ is degenerated to
$\tilde F_{1,(a,b)}$ -- a union of planes in the following configuration.

 \begin{center}
\epsfig{file=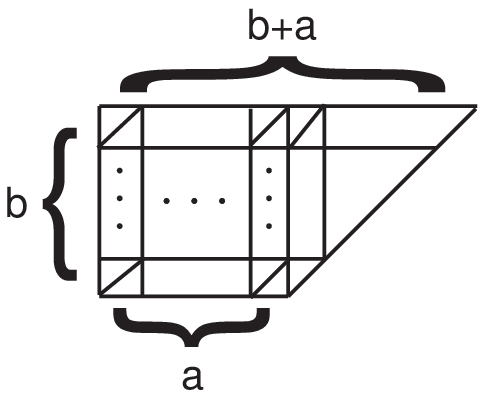}\\
(figure 1.5)
\end{center}

Each triangle represents a plane; each inner edge represents an
intersection line between planes.  The number of the planes is
$2ab+b^2$; the number of intersection lines is
$3ab-a+\frac{3b}{2}(b-1).$ We take a generic projection of $\tilde
F_{1,(a,b)}$ onto $\C\P^2$ where each plane is projected onto
$\C\P^2.$ The ramification curve of this projection is the union of
lines. The singular points of the ramification curve are represented
by the vertices. The branch curve of $\tilde F_{1,(a,b)}\to\C\P^2$,
denoted by $\tilde S_{(a,b)}$, is the image of the union of lines
and its singular points are the images of the vertices and the
intersection points in $\C\P^2$ of the images of any two of the
intersection lines. Special notations of the vertices and the edges
of the complex in Fig.~1.5 (which represent $\tilde S_{(a,b)}$) will
be given in Section 4.

\section{$B_n,\tilde B_n$ and $\tilde B_n$-groups}

The aim of this section is to introduce a few facts about $B_n$ and a certain
quotient of it, which will serve us in the next section.

\begin{defn}\label{defn3.1} $B_n,S_n$:

The braid group on $n$ strings is
$$B_n=\Big\{x_1,\dots,x_{n-1}\Big| \begin{array}{ll}[x_1,x_j]=1\quad |i-j|>1  \\
 \langle x_i,x_j\rangle=1\quad |i-j|-1\end{array}\Big\}.$$
Recall that the permutation group is
$$S_n=\Big\{x_1,\dots,x_{n-1}\Big| \begin{array}{ll}[x_1,x_j]=1\quad |i-j|>1  \\
 \langle x_i,x_j\rangle=1\quad |i-j|-1\end{array}, x_i^2=1\Big\}.$$
So, $\exists$ homomorphism $\varphi: B_n\to S_n$. Denote by $\delta$
the degree homomorphism $\delta: B_n\to\Z$; denote
$P_n=\ker\varphi,$ \ $P_{n,0}=P_n\cap \ker\delta.$
\end{defn}

We now recall another definition of $B_n$.

Let $D$ be a closed disk in $ \mathbb{R}^2,$ \ $K\subset Int(D),$
$K$ finite, $n= \#K$. Recall that the braid group $B_n[D,K]$ can be
defined as the group of all equivalent diffeomorphisms $\beta$ of
$D$ such that $\beta(K) = K\,,\, \beta |_{\partial D} =
\text{Id}\left|_{\partial D}\right.$. \\

\begin{defn}  $H(\sigma)$, half-twist defined by
$\sigma$

Let $a,b\in K,$ and let $\sigma$ be a smooth simple path in $Int(D)$
connecting $a$ with $b$ \ s.t. $\sigma\cap K=\{a,b\}.$ Choose a
small regular neighborhood $U$ of $\sigma$ contained in $Int(D),$
s.t. $U\cap K=\{a,b\}$. Denote by $H(\sigma)$ the diffeomorphism of
$D$ which switches $a$ and $b$ by a counterclockwise 180 degree
rotation and is the identity on $D\setminus U$\,. Thus it defines an
element of $B_n[D,K],$ called {\it the half-twist defined by
$\sigma$ }.

\end{defn}

\begin{defn}\label{defn3.2} $\tilde B_n$

Let $\tilde B_n$ be the quotient of $B_n$ by the following
commutator, $\tilde B_n=B_n/\langle[x_2,(x_2)_{x_1\,x_3}]\rangle$,
that is, by the commutator of two half-twists intersecting
transversally.
\end{defn}

\begin{lemma}\label{lem3.1} Let $x,y\in\tilde B_n.$
\begin{enumerate}\item[(i)] If the endpoints of $x$ and $y$ are disjoint, then
$[x,y]=1.$ \item[(ii)] If the endpoints of $x$ and $y$ have one common
endpoint, the $\langle x,y\rangle=1.$
\end{enumerate}
\end{lemma}

\begin{proof} \cite[Section 3]{BGT5}.
\end{proof}Let $\tilde \varphi$ be the induced
homomorphism from $\varphi,$ s.t. $\tilde\varphi:\tilde B_n\to S_n$.
Define $\tilde P_n=\ker\tilde \varphi,$\ $\tilde
P_{n,0}=\ker\tilde\varphi\cap\ker\tilde\delta$ (where
$\tilde\delta:\tilde B_n\to\Z).$

We cite now the main results of \cite[Section 1]{Mo}; see also
\cite{Te2}.

\begin{lemma}\label{lem3.2} Denote by $x_i$ the image of the generator $X_i$ in
$\tilde B_n.$ Let $s_1=x_1^2,$\ $\mu=[x_1^2,x_2^2],$\
$u_i=[x_i^{-1},x_{i+1}^2]$ \ $\forall\ 1\le i\le
n-2,u_{n-1}=[x_{n-2}^2,x_{n-1}].$ So $\tilde P_{n,0}$ is generated
by $u_1,\dots, u_{n-1},$ and $\tilde P_n$ is generated by
$s_1,u_1,\dots,u_{n-1}.$

We also have the following:
$$[u_i,u_j]=\begin{cases} 1\quad & |i-j|>1\\
\mu\quad & \text{otherwise}\end{cases}$$
$$[s_1,u_i]=\begin{cases} 1\quad & i\ne 2\\
\mu\quad & i=2\end{cases}$$ Moreover, $\mu^2=1,$\
$\mu\in\Center(\tilde B_n)$ and $\langle\mu\rangle=[\tilde
P_{n,0},\tilde P_{n,0}]=[\tilde P_n,\tilde P_n].$ Therefore, $\tilde
P_{n,0}$ is solvable and $Ab(\tilde P_n)\simeq \Z^n$,\ $Ab(\tilde
P_{n,0})\simeq \Z^{n-1}.$

We can also formulate the action of $\tilde B_n$ on $\tilde P_n$ by
conjugation:
$$(s_1)_{x_i}=\begin{cases} s_1\quad & i\ne 2\\
s_1u_2^{-1}\quad & i=2\end{cases},\qquad (u_j)_{x_i}=\begin{cases} u_j\quad &
|i-j|>1\\
u_iu_j\quad & |i-j|=1\\
u_i^{-1}\mu\quad & i=j\end{cases}$$
\end{lemma}

Actually, this action of $\tilde B_n$ on $\tilde P_n$ was developed
(see \cite{BGT5}) to abstract groups with $\tilde B_n$ actions
similar to the action on $\tilde P_n$ and $P_{n,0}.$ This is
explained in the following properties.

\begin{defn}\label{defn3.3} \textit{Adjacent half-twists}

If $x,y\in\tilde B_n$ are two half-twists whose endpoints have only
one point in common (and they can intersect each other
transversally), We say $x$ and $y$ are adjacent.
\end{defn}

The following definitions, lemmas and propositions are taken from
\cite{Te2}.

\begin{defn}\label{defn3.4} \textit{Polarized half-twists, polarization}

We say that a half-twist $X\in B_n$ (or $\tilde X$ in $\tilde B_n$) is
polarized if we choose an order on the endpoints of $X.$ The order is called
the polarization of $X$ or $\tilde X.$
\end{defn}

\begin{defn}\label{defn3.5} \textit{Orderly adjacent}

Let $X,Y$ be two adjacent polarized half-twists in $B_n$ (resp. in $\tilde
B_n).$ We say that $X,Y$ are \textit{orderly adjacent} if their common point is
the ``end" of one of them and the ``origin" of another.
\end{defn}

The following definition derives its motivation from the action of $\tilde B_n$
on $\tilde P_n.$

\begin{defn}\label{defn3.6} $\tilde B_n$-\textit{group}

A group $G$ is called a $\tilde B_n$-group if there exists a homomorphism
$\tilde B_n\to \Aut(G).$ We denote $(g)_b$ by $g_b.$
\end{defn}

\begin{defn}\label{defn3.7} \textit{Prime element, supporting half-twist (s.h.t.)
corresponding central element}

Let $G$ be a $\tilde B_n$-group.

An element $g\in G$ is called a \textit{prime element} of $G$ if
there exists a half-twist $X\in B_n$ and $\tau\in\Center(G)$ with
$\tau^2=1$ and $\tau_b=\tau\ \forall \ b\in\tilde B_n$ such that
\begin{enumerate}
\item[(1)] $g_{\tilde X^{-1}}=g^{-1}\tau$ \item[(2)] For every half-twist $Y$
adjacent to $X$ we have:\\$g_{\tilde X \tilde Y^{-1}\tilde X^{-1}}=g_{\tilde
X}^{-1}g_{\tilde X\tilde Y^{-1}}$ \item[] $g_{\tilde Y^{-1}\tilde
X^{-1}}=g^{-1}g_{\tilde Y^{-1}}.$ \item[(3)] For every half-twist $Z$ disjoint
from $X,$\ $g_{\tilde Z}=g.$
\end{enumerate}

The half-twist $X$ (or $\tilde X$) is called the \textit{supporting half-twist}
of $g$\ ($X$ is the s.h.t. of $g.$)

The element $\tau$ is called the \textit{corresponding central element}.
\end{defn}

\begin{lem}\label{lem3.3}
Let $G$ be a $\tilde B_n$-group.

Let $g$ be a prime element in $G$ with supporting half-twist $X$ and
corresponding central element $\tau.$ Then:
\begin{enumerate}
\item[(1)] $g_{\tilde X}=g_{\tilde X^{-1}}=g^{-1}\tau,$\ $g_{\tilde X^2}=g.$
\item[(2)] $g_{\tilde Y^{-2}}=g\tau\ \forall\ Y$ consecutive half-twist to $X.$
\item[(3)] $[g,g_{\tilde Y^{-1}}]=\tau\ \forall\ Y$ consecutive half-twist to
$X.$\end{enumerate}
\end{lem}

\begin{defn}\label{defn3.8} \textit{Polarized pair}

Let $G$ be a $\tilde B_n$-group, \ $h$ a prime element of $G,$\ $X$ its
supporting half-twist.  If $X$ is polarized, we say that $(h,X)$ (or $(h,\tilde
X)$) is a polarized pair with central element $\tau,$\ $\tau=hh_{\tilde
X^{-1}}.$
\end{defn}

\begin{defn}\label{defn3.9} \textit{Coherent pairs, anti-coherent pairs}

We say that two polarized pairs $(h_1,\tilde X_1)$ and $(h_2,\tilde X_2)$ are
coherent (anti-coherent) if $\exists \tilde b\in \tilde B_n$ such that
$(h_1)_{\tilde b}=h_2,$\ $(\tilde X_1)_{\tilde b}=\tilde X_2,$ and $\tilde b$
preserves (reverses) the polarization.
\end{defn}

\begin{prop}\label{prop3.1} Let $(h,\tilde X)$ be a polarized pair, $h\in G,$\
$\tilde X\in\tilde B_n.$ Let $\tilde T$ be a polarized half-twist in $\tilde
B_n.$ Then there exists a unique prime element $g\in G$ such that $(g,\tilde
T)$ and $(h,\tilde X)$ are coherent.
\end{prop}

\begin{defn}\label{defn3.10} $L_{h,\tilde X}(\tilde T)$

Let  $(h,X)$ be a polarized pair $\tilde T\in\tilde B_n.$ \ $L_{h,\tilde
X}(\tilde T)$ is the unique prime element s.t. $(L_{(h,\tilde X)}(\tilde
T),\tilde T)$ is coherent with $(h,\tilde X).$
\end{defn}

In fact, one can prove that $\tilde P_n$ (as a $\tilde B_n$-group) has a prime
element, and $\tilde P_{n,0}$ is generated by the orbit of this prime element.

\begin{lem}\label{lem3.4} Let $X_1,X_2$ be $2$ consecutive half-twists in
$B_n.$ Let $u=(\tilde X_1^2)_{\tilde X_2^{-1}}\tilde X_2^{-2}$. Then
$u\in\tilde P_{n,0},$\ $u$ is  a prime element in $\tilde P_n$ (considered as a
$\tilde B_n$-group), and $\tilde X_1$ is the supporting half-twist of $u.$
\end{lem}

\begin{lem}\label{lem3.5}
$\tilde P_{n,0}$ is  a primitive $\tilde B_n$-group generated by the $\tilde
B_n$-orbit of a prime element $u=\tilde X^2\tilde Y^{-2}$, where $\tilde
X,\tilde Y$ are adjacent half-twists in $\tilde B_n,$\ $\tilde T=\tilde X\tilde
Y\tilde X^{-1}$ is a supporting half-twist for $u.$
\end{lem}

We shall also cite from \cite{Te2} the criterion for prime elements
in $\tilde B_n$-groups; we will not use it directly, but rather
implicitly, when quoting, in Section 4, the results for the $\tilde
B_n$-groups (see Lemma 4.2).

\begin{prop}\label{prop3.2} Assume $n\ge5.$ Let $G$ be a $\tilde B_n$-group,
and let
$$(\tilde X_1,\tilde X_2,\dots,\tilde X_{n-1})$$
be a standard base of $\tilde B_n.$ Let $S$ be an element of $G$ with the
following properties:
\begin{enumerate}
\item[(0)] $G$ is generated by $\{S_b,\ b\in\tilde B_n\};$ \item[(1$_a$)]
$S_{\tilde X_2^{-1}\tilde X_1^{-1}}=S^{-1}S_{\tilde X_2^{-1}};$
 \item[(1$_b$)]
$S_{\tilde X_1\tilde X_2^{-1}\tilde X_1^{-1}}=S^{-1}_{\tilde X_1}S_{\tilde
X_1,\tilde X_2^{-1}};$ \item[(2)] For $\tau=SS_{\tilde X_1^{-1}},$\
$T=S_{\tilde X_2^{-1}},$ we have:
\begin{itemize} \item[(2$_a$)] $\tau_{\tilde X_1^2}=\tau;$
\item[(2$_b$)] $\tau_T=\tau_{\tilde X_1}^{-1};$
\end{itemize}
\item[(3)] $S_{\tilde X_j}=S\ \forall j\ge3;$ \item[(4)] $ S_c=S,$ where
$c=[\tilde X_1^2,\tilde X_2^2].$
\end{enumerate}
Then $S$ is a prime element of $G,$\ $\tilde X_1$ is a supporting half-twist of
$S$ and $\tau$ is the corresponding central element. In particular,
$\tau^2=1,$\ $\tau\in\Center(G),$\ $\tau_b=\tau\ \forall\ b\in\tilde B_n.$
\end{prop}

\section{Calculation of the fundamental group}

In this section we will calculate the fundamental group of the
complement of the branch curve of $F_{1,(a,b)}$. This  computation
requires explicit knowledge of the braid monodromy factorization
(BMF) technique. This knowledge can be found at
\cite{KuTe},\cite{MT2}, \cite{MT3}. However, we recall the main
definitions regarding the braid monodomy factorization related to a
curve $S$.

 \begin{defn} The braid monodromy w.r.t. $S,\pi,u$

 Let $S$ be a curve, $S\subseteq \C^2$
 Let $\pi: S\to\C^1$ be defined by
 $\pi(x,y)=x.$
We denote $\deg\pi$ by $m.$ Let $N=\{x\in\C^1\bigm| \#\pi^{-1}(x)<
m\}.$
   Take $u\notin N,\, u$ real, s.t.  $\Re(x)\ll u$ \ $\forall x\in N.$
Let  $ \C^1_u=\{(u,y)\}.$  There is a  natural defined homomorphism
$$\pi_1(\C^1-N,u)\xrightarrow{\vp} B_m[\C_u^1,\C_u^1\cap S]$$ which
is called {\it the braid monodromy w.r.t.} $S,\pi,u,$ where $B_m$ is
the braid group. We sometimes denote $\vp$ by $\vp_u$. Note that in
this definition we regard $B_m$ as the group of diffeomorphisms, as
described in the
previous section.\\

Denote the generator of the center of $B_n$ as $\Delta^2$. We recall
Artin's theorem on the presentation of $\Delta^2$ as a product of
braid monodromy elements of a geometric-base (a base of $\p =
\p(\C^1 - N, u)$ with certain
properties; see \cite{BGT1} for definitions).\\
 \textbf{Theorem}: Let $S$ be a curve transversal to
the line in infinity, and $\vp$ is a braid monodromy of $S , \vp:\p
\rightarrow B_m$. Let {$\delta_i$} be a geometric (free) base
(g-base) of $\p.$ Then:
$$\Delta^2 = \prod\vp(\delta_i).$$ This product is also defined as
the \textsl{braid monodromy factorization} (BMF) related to a curve $S$.\\

\end{defn}

Since $\tilde S_{F_{1,(a,b)}}$, which is the branch curve of the
degenerated surface $\tilde F_{1,(a,b)}$, is a line arrangement, we
can compute the braid monodromy factorization as in \cite{BGT1}. In
order to compute the braid monodromy factorization of
$S_{F_{1,(a,b)}}$, we use the regeneration rules (\cite{MT3}).
 The regeneration methods are actually, locally, the
reverse process of the degeneration method. When regenerating a
singular configuration consisting of lines and conics, the final
stage in the regeneration process involves doubling each line, so
that each point of $K$ (which is the set of points in the disk, that
is $\C_u^1\cap \tilde S_{F_{1,(a,b)}}$) corresponding to a line
labelled $i$ is replaced by a pair of points, labelled $i$ and $i'$.
The purpose of the regeneration rules is to explain how the braid
monodromy behaves when lines are doubled in this manner.\\

Let $F_{1,(a,b)}$, \ $a,b>1$ be the Hirzebruch surface embedded
w.r.t. the linear system\\ $|aC+bE_0|.$ As shown, $F_{1,(a,b)}$
could be degenerated into a union of $2ab+b^2$ planes in the
following arrangement:

 \begin{center}
\epsfig{file=./fig1_5.EPS}\\
(figure 4.1)
\end{center}
We shall give a special presentation of $B_n,$ from which we will
induce an injection of $\tilde B_n$ to
$G=\pi_1(\C^2-S_{F_{1,(a,b)}}).$

 \textbf{Remark}: From now on, we
denote by $\bar{S}_{F_{1,(a,b)}}$  the branch curve of $F_{1,(a,b)}$
(in $\C\P^2$), and by $S_{F_{1,(a,b)}}$ a generic affine portion of
it (in $\C^2$).

 Let $a,b$ be integers $b>1, a\geq1$\ $n=2ab+b^2$.
Let $s_{ij}=(i,j),$\ $t_{ij}=\left(i+\frac{1}{2},j\right)\in\mathbb
R^2.$ Let $K_{a,b}$ be the set in $\mathbb R^2$ consisting of the
points $s_{ij},t_{ij},$\ $1\le j\le b,$\ $1\le i\le a+j;$ so $\#
K_{a,b}=2ab+b^2.$

Let $D$ be a large disk in $\mathbb R^2$ containing $K_{a,b}.$
Consider the oriented line segments $\vec x_{ij}=[s_{i,j},t_{i,j}]$
where $1\le j\le b,$ \ $1\le i\le a+j;$\ $\vec
y_{ij}=[t_{i,j},s_{i+1,j}],$\ $1\le j\le b,$\ $1\le i\le a+j-1;$\
$\vec z_{ij}=[s_{i,j},t_{i,j+1}],$\ $1\le j\le b-1,$\ $1\le i\le
a+j.$ Consider $B_n=B_n[D,K_{a,b}].$ Let $X_{ij},Y_{ij},\underline
Z_{ij}$ be polarized half-twists in $B_n$ by the oriented segments
$\vec x_{ij},\vec y_{ij},\vec z_{ij}$ respectively. Let
$Z_{ij}=\underline Z_{ij}$, when $i=a+j,$\ $1\le j\le b-1.$ We
define $Z_{ij}$ for $1\le i\le a+j-1,$ inductively:
$$Z_{ij}=(Z_{i+1,j})_{X_{i+1,j+1}^{-1}Y_{i,j+1}Y_{i,j}^{-1}X_{i,j}}.$$

\begin{prop}\label{prop4.1} $B_n$ can be finitely presented as follows:

 \noindent Generators:

\parindent 60pt $X_{i,j},\ 1\le j\le b,\ 1\le i\le a+j.$

 $Y_{ij},$\ $1\le j\le b; \ 1\le i\le a+j-1$.

  $Z_{ij},$\
$i=a+j,\ 1\le j\le b-1.$

\noindent Relations:

\parindent 60pt $\forall$\ two generators $a,b$ of the above which are
adjacent, $\langle a,b\rangle=1.$

$\forall$\ two generators $c,d$ which are disjoint $[c,d]=1.$

$\forall \ j\in(1,\dots,b-1),$\ $i=a+j:[X_{i,j},Z_{ij}Y_{i-1,j}Z_{ij}^{-1}]=1.$
\end{prop}

\begin{proof} This is a standard consequence of the usual presentation of
$B_n[D,K_{a,b}]$ (see \cite{BGT1}).

\end{proof}

The formulas define inductively a polarization for each $Z_{ij}.$
One can check that it coincides with the given polarization of
$\underline Z_{ij}$, i.e., corresponds to the ordered pair
$(s_{ij},t_{i,j+1}).$

Denote by $x_{ij},y_{ij},z_{ij}$ the images of
$X_{ij},Y_{ij},Z_{ij}$ in $\tilde B_n.$ Thus we get a representation
of $\tilde B_n.$ We consider $\{x_{ij},y_{ij},z_{ij}\}$ with
polarization introduced above.

\begin{defn}\label{defn4.1}
Let $G$ be a primitive $\tilde B_n$-group generated by the orbit of
a prime element $B_{1,1}$ supported by the half-twist $Y_{1,1}.$
According to Proposition 3.1, $\forall$ polarized half-twist
$t\in\tilde B_n$,\ $\exists$\ unique prime element
$L_{\{B_{1,1},y_{1,1}\}}(t)\in G,$ s.t. the pair
$\{L_{\{B_{1,1},y_{1,1}\}}(t),t\}$ is coherent with
$\{B_{1,1},y_{1,1}\}.$

 Define
$$A_{ij}=L_{\{B_{1,1},y_{1,1}\}}(x_{ij})$$
$$B_{ij}=L_{\{B_{1,1},y_{1,1}\}}(y_{ij})$$
$$C_{ij}=L_{\{B_{1,1},y_{1,1}\}}(z_{ij})$$ \end{defn}

\begin{remark}\label{rem4.1}
Looking at \cite[Remark 6]{Mo}, one gets the formulas for the $\tilde
B_n$-action on $G$ in terms of $\{x_{ij},y_{ij},z_{ij};\ A_{ij},B_{ij},C_{ij};\
i,j=\dots\}.$ In particular, we see that $G$ is generated by
$\{A_{ij},B_{ij},C_{ij}\}$ (because $G$ is generated by the $\tilde B_n$-orbit
of $B_{1,1}).$

\end{remark}

Denote by $\tilde S_{a,b}:=\tilde S_{F_{1,(a,b)}}$ the degenerated branch curve
of $F_{1,(a,b)}$. We define now a planar 2-complex, to represent the polygon in
Fig.~4.1.

\begin{defn}\label{defn4.2} We use a planar 2-complex $K(a,b)$ defined as
follows: $K(a,b)\subset\mathbb R^2.$ Define $P$, the polygon whose vertices are
$(0,0),(a,0),(a+b,b),(0,b).$ So the vertices of $K(a,b)$ are the points
$\left\{\omega_{rk}=(r,k)\bigm|_{r,k\in\mathbb Z}^{\omega_{rk\in P}}\right\}.$
The edges of $K(a,b)$ are the straight line segments of the following three
types:
\begin{enumerate}\item[(a)] ``diagonal":\quad $[\omega_{r,k},\omega_{r+1,k+1}],$\
$0\le k\le b-1,$\ $0\le r\le a+k;$ \item[(b)] ``vertical":\quad\
$[\omega_{r,k},\omega_{r,k+1}],$\ $0\le k\le b-1,$\ $0\le r\le a+k;$
\item[(c)] ``horizontal":\,$[\omega_{r,k},\omega_{r+1,k}],$\ $0\le
k\le b,$\ $0\le r\le a+k-1;$
\end{enumerate}
The 2-simplices of $K(a,b)$ are the triangles
$\Delta\{\omega_{r,k},\omega_{r+1,k},\omega_{r+1,k+1}\}$ and
$\Delta\{\omega_{rk},\omega_{r,k+1},\omega_{r+1,k+1}\}.$
\end{defn}

\begin{defn}\label{defn4.3}
The vertices $\omega_{rk}$ that are not on the boundary of $P$ will
be called 6-point; the vertices $\omega_{0,0},\omega_{a,0}$ will be
called 2-point; and all the other vertices $\omega_{rk}$ on the
boundary of $P$ s.t. $(r,k)\ne(0,b),(a+b,b)$ will be called 3-point.
\end{defn}

\begin{defn}\label{defn4.4}

\textbf{\textsl{(1)}} Consider $B_m=B_m[D,K]$, where $D$ is a large
disk in $\C^1,$ centered\\ at (0) and
\begin{alignat*}{2}
&K=&&\{q_{rk\delta}^{(\varepsilon)}\bigm|\varepsilon=1,2,3,\ \delta=0,1\
\text{s.t.}:\\
  & &&\text{for}\ \varepsilon=1,\ 1\le k\le b,\ 1\le r\le a+k-1\\
  & &&\text{for}\ \varepsilon=2,\ 1\le k\le b,\ 1\le r\le a+k-1\\
  & &&\text{for}\ \varepsilon=3,\ 1\le k\le b-1,\ 1\le r\le a+k\\
  & && q_{rk\delta}^{(\varepsilon)}\ \text{are real points such that}\
 q_{rk0}^{(\varepsilon)},q_{rk1}^{(\varepsilon)} \ \text{are very close to each other,
  and}\\
   & && q_{rk\delta}^{(\varepsilon)}<q_{r'k'\delta'}^{(\varepsilon')}\ \text{if
    either}\
  k<k'\ \text{or}\ k=k'\ \text{and}\ r<r'\\& &&\text{or}\ k=k', \  r=r'\ \text{and}\
  \varepsilon<\varepsilon',
  \ \text{or}\ k=k',  \ r=r', \ \varepsilon=\varepsilon'\ \text{and } \delta<\delta'
  \}.
 \end{alignat*}

 The points of $K$ that we associate with the non-boundary edges of $K(a,b)$ are
 as follows: $q_{rk0}^{(1)},q_{rk1}^{(1)}$ correspond to the diagonal edge
 $[\omega_{r-1,k-1},\omega_{r,k}];$\ $q_{rk0}^{(2)},q_{rk1}^{(2)}$ correspond
 to the vertical edge $[\omega_{r,k-1},\omega_{r,k}];$ \
 $q_{rk0}^{(3)},q_{rk1}^{(3)}$ correspond to the horizontal edge
 $[\omega_{r-1,k},\omega_{r,k}].$

 As was indicated earlier, during the regeneration process, each
 line doubles itself, and thus each point of $\C^1\cap \tilde S_{a,b}$
 is replaced by a pair of points, which are
 $q^{(\varepsilon)}_{rk0}$ and $q^{(\varepsilon)}_{rk1}$.

 \textbf{\textsl{(2)}} Let
 $$m_{r,k}=\begin{cases} 12\quad & \text{if}\quad \omega_{rk}\ \text{is a}\
 6\text{-point}\\
 4\quad & \text{if}\quad \omega_{rk}\ \text{is a}\
 3\text{-point}\\
 2\quad & \text{if}\quad \omega_{rk}\ \text{is a}\
 2\text{-point}
 \end{cases}$$

 Denote by $K_{r,k}$ the subset of $K$ consisting of the points associated with
 the non-boundary edges of $K(a,b)$ which meet at $\omega_{r,k}.$ Clearly,
 $\#K_{r,k}=m_{r,k}.$

 \textbf{\textsl{(3)}} Denote $f_{rk}: B_{m_{r,k}}\to B_m[D,K]$ an embedding of $B_{m_{r,k}}$
 into $B_m[D,K]$ corresponding to a connection below the real axis of the points of $K_{r,k}$ by consecutive
 simple paths (see \cite{BGT1}). Clearly, each $B_{m_{r,k}}$ is either
 $B_{12},B_4$ or $B_2.$

  \end{defn}

 From each 6/3/2-point, relations between the generators of the fundamental
 group $\pi_1(\C^2-S_{F_{1,(a,b)}})$ can be induced. These relations are
 written with the same notations as in \cite{Mo}. We refer the reader to this
 article.  However, we state a few of the main results.

 Consider $K\subset D,$\ $K=\{q_{rk\delta}^{(\varepsilon)}\}.$ Take a point
 $\underline u$ on $\partial D$ below the real axis. Using small (positively
 oriented) circles around the points $q_{rk\delta}^{(\varepsilon)}$ and
 connecting these circles by (straight) simple lines with $\underline u,$ we
 obtain a geometric base $\{\gamma_{rk\delta}^{(\varepsilon)}\}$ for
 $\pi_1(D-K,\underline u).$

 A full set of relations between $\{\gamma_{rk\delta}^{(\varepsilon)}\}$ can be
 described, corresponding to the braid monodromy factorization (see \cite{BGT1} for a formula
 computing the BMF of a generic line arrangement - which is actually the factorization
 on which we perform the regeneration process to get the following):
 $$\Delta^2 = \varepsilon(a,b)=\prod_{\omega_{r,k}}C(r,k)\mathcal{H}(r,k),$$
 where $\mathcal{H}(r,k)$ are the factorizations induced from the
 6/3/2-points -
 $\omega_{r,k}\,$(see Appendix).$\,\,C(r,k)$ are the factorizations that we get from the parasitic
 intersection of the branch curves (see \cite[Chapter~2]{Mo} or \cite{BGT1}). We get
 a presentation of $\pi_1(\C^2-S_{F_{1,(a,b)}})$ by using the Van-Kampen
 Theorem \cite{vK} which says that from each factor from $\varepsilon(a,b)$, a
 relation between $\{\gamma_{rk\delta}^{(\varepsilon)}\}$ can be induced.
 Taking a
braid which is a half-twist that corresponds to a path $\sigma$ from
$q_{r_1k_1\delta_1}^{\varepsilon_1}$ to
$q_{r_2k_2\delta_2}^{\varepsilon_2}$ via $u,$ we let $\delta_1$
(resp. $\delta_2$) be the path from $u$ to
$q_{r_1k_1\delta_1}^{\varepsilon_1}$ (resp.
$q_{r_2k_2\delta_2}^{\varepsilon_2}$) along $\sigma$, going around
$q_{r_1k_1\delta_1}^{\varepsilon_1}$ (resp.
$q_{r_2k_2\delta_2}^{\varepsilon_2}$) and coming back to $u$ along
the same path, respectively. Let $A$ and $B$ be the homotopy classes
of a loop around $q_{r_1k_1\delta_1}^{\varepsilon_1}$ (resp.
$q_{r_2k_2\delta_2}^{\varepsilon_2}$) along $\delta_1$ (resp.
$\delta_2$). $A$ (resp. $B$) is a conjugation of
$\gamma_{r_1k_1\delta_1}^{\varepsilon_1}$ (resp.
$\gamma_{r_2k_2\delta_2}^{\varepsilon_2}).$

By the Van Kampen Theorem, we have one of the following relations in
$\pi_1(\C^2-S_{F_{1,(a,b)}})$ (fixed according to the type of
singularity, from which we have the path $\sigma$):
\begin{enumerate}
\item[1.] $A=B,$ if the singularity is a branch point, \item[2.]
$[A,B]=ABA^{-1}B^{-1}=e$ if it is a node, \item[3.] $\langle
A,B\rangle=ABAB^{-1}A^{-1}B^{-1}=e$ if it is a cusp.
\end{enumerate}

 \begin{defn}\label{defn5.5} Let
 \begin{align*}
 &\ell_{r,k}^{(1)}=\begin{cases} 1-k\quad & \text{for}\quad r\ge k\\
 1-r\quad &\text{for}\quad r<k
\end{cases}\\
& \ell_{r,k}^{(2)}=k-1\\
&\ell_{r,k}^{(3)}=0.
\end{align*}
(Evidently, $\ell_{r+1,k}^{(3)}=\ell_{r,k}^{(3)};$\
$\ell_{r,k+1}^{(2)}=\ell_{r,k}^{(2)}+1;$\
$\ell_{r+1,k+1}^{(1)}=\ell_{r,k}^{(1)}-1.)$ Let
$e_{rk\delta}^{(\varepsilon)}=(\gamma_{rk\delta}^{(\varepsilon)})(\rho_{rk}^{(\varepsilon)})^{\ell_{rk}^{(
\varepsilon)}}$ (where $\rho_{rk}^{(\varepsilon)}$ is the half-twist
in $B_m[D,K]$ defined by the segment
$[q_{rk0}^{(\varepsilon)},q_{rk1}^{(\varepsilon)}].$
\end{defn}

\begin{defn}\label{defn4.6}
Denote by $G$ the group defined by $\varepsilon(a,b);$ more
precisely, the quotient of the free group generated by
$\{e_{rk\delta}^{(\varepsilon)}\},$ modulo relations (we call them
$R\varepsilon)$ induced from 6/3/2-points, and the relation induces
from the parasitic intersections, for all $\omega_{r,k}$ (see
\cite[Chapter~3]{Mo} for those relations or in the Appendix).

By the definition of $\varepsilon(a,b)$ (braid monodromy
factorization for $S_{F_{1,(a,b)}}),$ we have\\ $G\simeq
\pi_1(\C^2-S_{F_{1,(a,b)}},\underline u).$ Let
$E_{rk\delta}^{(\varepsilon)}$ be the images of
$e_{rk\delta}^{(\varepsilon)}$ in $G.$
\end{defn}

\begin{prop}\label{prop4.2} $\exists$ homomorphism $\tilde\alpha:\tilde B_n\to
G$ which is defined by:
$$\tilde \alpha(x_{ij})=E_{ij0}^{(1)} ,\quad \tilde
\alpha(y_{ij})=E_{ij0}^{(2)}\quad \forall i,j,$$$$\quad \tilde
\alpha(z_{ij})=E_{ij0}^{(3)}\quad\text{(where}\ i=a+j);$$ moreover,
$$\tilde\alpha(z_{ij})=E_{ij0}^{(3)}\ \forall\ (i,j)\in\Vertices(K(a,b)),\ i\ne a+j.$$
\end{prop}

\begin{proof} See \cite[Proposition~8]{Mo}. See the induced relations for each
2/3/6-point and explanation why $\tilde B_n$ can be embedded in $G$
in the Appendix.
\end{proof}

Let $E_{rk}^{(\varepsilon)}=E_{rk0}^{(\varepsilon)}$,\ $\mathcal{B}$
be the subgroup of $G$ generated by $\{E_{rk}^{(\varepsilon)}\}.$ It
follows from Proposition~\ref{prop4.2} that $\B=\tilde\alpha(\tilde
B_n).$ Let $\PP=\tilde\alpha(\tilde P_n),$\
$\PP_0=\tilde\alpha(\tilde P_{n,0})$ (where $P_{n,0}=\ker(P_n\to
Ab(B_n)),$\ $\tilde P_{n,0}$ is the image of $P_{n,0}$ in $\tilde
B_n).$ From \cite[Theorem~1]{Mo} or from Lemma 3.3,  it follows that
$\tilde P_{n,0}$ is a primitive $\tilde B_n$-group with prime
element $u=(y_{1,1}^2)_{x_{2,1}^{-1}}x_{2,1}^{-2}$ ($x_{2,1}$ and
$(y_{1,1})_{x_{2,1}^{-1}}$ are two adjacent half-twists in $\tilde
B_n$), and s.h.t. equal to $y_{11}.$ Denote
$c=[y_{1,1}^2,x_{2,1}^2]$. We get from \cite[Theorem ~1]{Mo} that
$c^2=1,$ \ $c\in\Center(\tilde B_n)$ and $c$ generates $\tilde P_n'$
and $\tilde P_{n,0}'$.  Denoting $\eta_{1,1}=\tilde
\alpha(u)=(E_{1,1}^{(2)})_{(E_{2,1}^{(1)})^{-1}}^2\cdot(E_{2,1}^{(1)})^{-2}$,\
$\mu=\tilde \alpha(c)=[(E_{1,1}^{(2)})^2,(E_{2,1}^{(1)})^2],$ we get
that $\PP_0$ is a primitive $\tilde B_n$-group, $\eta_{1,1}$ is a
prime element of $\PP_0$ with s.h.t. $y_{1,1},$ $\mu^2=1$,
$\mu\in\Center(\B)$ and $\mu$ generates $\PP'$ and $\PP_0'.$ Using
the polarization of $X_{i,j} Y_{i,j}, Z_{i,j}$ and Proposition
\ref{prop4.1}, we can find $\forall\ t\in\{x_{ij},y_{ij},z_{ij}\}$
(the generators of $\tilde B_n)$ and $\{z_{i,j}\bigm|
(i,j)\in\Vertices(K(a,b)),\ i,j\ge 1,\ i\ne a+j\}$ unique
$L_{\{\eta_{1,1},y_{1,1}\}}(t)\in \PP_0$ s.t. the pair
$\{L_{\{\eta_{1,1},y_{1,1}\}}(t),t\}$ is coherent with
$\{\eta_{1,1},y_{1,1}\}.$

\begin{defn}\label{defn4.7} Recall that $u=(y_{11}^2)_{x_{2,1}^{-1}}x_{2,1}^{-2},$ \
$\eta_{1,1}=\tilde\alpha(u).$ Define
$$\xi_{i,j}=L_{\{\eta_{1,1},y_{1,1}\}}(x_{ij}),\quad
\eta_{i,j}=L_{\{\eta_{1,1},y_{1,1}\}}(y_{i,j})$$
$$\zeta_{i,j}=L_{\{\eta_{1,1},y_{1,1}\}}(z_{i,j}).$$
\end{defn}

\begin{lemma}\label{lem4.1} $\mu\in\Center(G).$
\end{lemma}

\begin{proof} See \cite[Lemma 16]{Mo}.
\end{proof}

\begin{defn}\label{defn4.8} Let
$$d_{rk}=E_{rk1}^{(1)}(E_{rk0}^{(1)})^{-1},\quad v_{rk}=E_{rk1}^{(2)}(E_{rk0}^{(2)})^{-1},\quad
h_{rk}=E_{rk1}^{(3)}(E_{rk0}^{(3)})^{-1}.$$ ($d,v,h$ correspond to
``diagonal", ``vertical", ``horizontal".) Clearly, $G$ is generated
by $\{d_{rk},v_{rk},h_{rk};\,r,k=\dots\}$ and $\B.$ Denote by
$\mathcal{H}$ the subgroup of $G$ generated by the $\B$- (or $\tilde
B_n$-) orbit of $v_{1,1}.$
\end{defn}

\begin{lemma}\label{lem4.2}
{}\quad

\begin{enumerate}\item[1)] $\mathcal{H}$ is a primitive $\tilde B_n$-group with prime
element $v_{1,1},$ s.h.t. $y_{1,1}.$ \item [2)] $v_{1,1}$ is
actually a prime element of $G$ with s.h.t. $y_{1,1}$ (i.e.,
$v_{1,1}\cdot(v_{1,1})_{y_{1,1}^{-1}}\in\Center(G)).$
\end{enumerate}
\end{lemma}

\begin{proof} As in \cite[Lemma 17]{Mo}.
\end{proof}

\begin{defn}\label{defn4.9}
Using the polarization of $X_{i,j},Y_{i,j},Z_{i,j}$, we find
$\forall \ t\in\{x_{ij},y_{ij},z_{ij}\}\exists\
\,!\,L_{\{v_{1,1},y_{1,1}\}}(t)\in \mathcal{H}$ s.t. the pair
$\{L_{\{v_{1,1},y_{1,1}\}}(t),t\}$ is coherent with $\{v_{1,1},
\,y_{1,1}\}.$ Define
$$a_{i,j}=L_{\{v_{1,1},y_{1,1}\}}(x_{ij}),\quad
b_{i,j}=L_{\{v_{1,1},y_{1,1}\}}(y_{ij}),\quad
c_{i,j}=L_{\{v_{1,1},y_{1,1}\}}(z_{ij}).$$
\end{defn}

\begin{remark}\label{rem4.2} $\xi_{ij}$, $\eta_{ij},$\ $\zeta_{ij}$
$(a_{ij},b_{ij},c_{ij})$ coincide with $A_{ij},$\ $B_{ij},$\
$C_{ij}$ introduced in Definition \ref{defn4.1} for an arbitrary
primitive $\tilde B_n$-group $G,$ when this $G$ is replaced by
$\PP_0$ (resp. $\mathcal{H}$), and $\{B_{11},Y_{11}\}$ is replaced
by $\{\eta_{11},y_{11}\}$ (resp. $\{v_{11},y_{11}\}).$ Therefore,
replacing $A_{ij},$\ $B_{ij}$\ $C_{ij}$ by $\xi_{ij}$,\
$\eta_{ij}$,\ $\zeta_{ij}$ (resp. $a_{ij},$\ $b_{ij}$, \ $c_{ij})$,
we obtain formulas expressing the $\tilde B_n$-action on $\PP_0$
(resp. on $\mathcal{H}$). In particular, $\PP_0$ (resp.
$\mathcal{H}$) is generated by $\{\xi_{ij},\eta_{ij},\zeta_{ij}\}$
(resp. $\{a_{ij},b_{ij},c_{ij}\}).$
\end{remark}

\begin{defn}\label{defn4.10} $\forall\ x_{i,j},y_{i,j},z_{i,j},$ let $\tilde
x_{i,j}=\tilde \alpha(x_{i,j}),\tilde y_{i,j}=\tilde\alpha(y_{i,j}),$\ $\tilde
z_{i,j}=\tilde\alpha(z_{i,j}).$
\end{defn}

\begin{remark}\label{rem4.3} We have, by \cite[Remark 30]{Mo}, the following:
\begin{align*}
&d_{r+1,k+1} =(d_{rk})_{\tilde z_{rk}\tilde y_{rk}\tilde z_{r+1,k}^{-1}\tilde
y_{r,k+1}^{-1}}\\
 &h_{r+1,k}=(h_{rk})_{\tilde x_{rk}^{-1}\tilde y_{rk}\tilde y_{r,k+1}^{-1}\tilde
 x_{r+1,k+1}}\\
 &v_{r,k+1}=(v_{rk})_{\tilde x_{rk}^{-1}\tilde z_{rk}\tilde z_{r+1,k}^{-1}\tilde
 x_{r+1,k+1}}\\
 &h_{rk}=(v_{rk}d_{rk}(v_{rk}^{-1})_{x_{rk}^{-1}})_{z_{rk}x_{rk}}\\
 &v_{rk}=(h_{rk}d_{rk}(h_{rk}^{-1})_{x_{rk}^{-1}})_{y_{rk}x_{rk}}.
\end{align*}
\end{remark}

\begin{remark}\label{rem4.4} By \cite[Remark 31]{Mo}, we have
\begin{align*}
&d_{r+1, 1} =  \tilde
y_{r1}^{-2}(v_{r,1})_{x_{r+1,1}^{-1}y_{r,1}^{-1}}\cdot(\tilde
y_{r,1}^2)_{x_{r+1,1}
^{-1}}\\
 &d_{1,k+1}=\tilde z_{1,k}^{-2}(h_{1,k})_{x_{1,k}^{-1}z_{1,k}^{-1}}\cdot(\tilde
 z_{1,k}^2)_{x_{1,k+1}^{-1}}\\
& v_{r,b}=\tilde
x_{r,k}^{-2}\cdot(d_{r,b})_{y_{r,b}^{-1}x_{r,b}^{-1}}\cdot(\tilde
 x_{r,b}^2)_{y_{r,b}^{-1}}
\end{align*}
\end{remark}

Notice that in the following calculation, we will use the fact that
$\mu^2=\nu^2=1$ (since they are central elements).

\begin{prop}\label{prop4.3} Let $\lambda(k)=\frac{k(k-1)}{2}.$ We have
\begin{alignat*}{2}
&h_{rk}=c_{rk}^k\zeta_{rk}^{-k+1}(\mu\nu)^{\lambda(k)}\quad &&\forall \ r,k\\
&d_{rk}=a_{rk}^{r-k}\xi_{rk}^{-r+k}(\mu\nu)^{\lambda(k-r)}\quad&&\forall\ r,k\\
&v_{rk}=b_{rk}^r\eta_{rk}^{-r+1}(\mu\nu)^{\lambda(r)}\quad&&\forall\ r<a
\end{alignat*}
\end{prop}

\begin{proof} See \cite[Proposition 10]{Mo}.
\end{proof}

\begin{prop}\label{prop4.4} $v_{a,k}=1,\, \forall 0 \leq k \leq b.$
\end{prop}

\begin{proof} By the definition,
 $v_{a,0}=E_{a01}^{(2)}(E_{a00}^{(2)})^{-1}$, but $\omega_{a,0}$ is a 2-point,
 and the induced relation from it is $\gamma_{a00}^{(2)}=\gamma_{a01}^{(2)}$ or
 $1=E_{a01}^{(2)}(E_{a00}^{(2)})^{-1}$. by the  relation $v_{r,k+1}=(v_{r,k})_{\tilde
 x_{rk}^{-1}z_{rk}z_{r+1,k}^{-1}\tilde x_{r+1,k+1}},$ we can see that
 $v_{a,k}=1$\ $\forall \ 0\le k\le b.$
 \end{proof}

\begin{prop}\label{prop4.4} For
$(r,k)\in\{(a+1,2),(a+2,3),\dots(a+b-1,b-1)\}=:I,$
$$v_{r,k}=(E_{r,k-1}^{(3)})^{-2}h_{r,k-1}^{-1}(h_{r,k-1})_{(E_{r,k}^{(2)})^{-1}}(E_{r,k-1}^{(3)})_{(E_{r,k}^{
(2)})^{-1}}^2.$$
\end{prop}

\begin{proof} Assume $(r,k)=(a+1,2).$ The proof for the other points is the
same.

We have by the relations induced from the 3-point $\omega_{a+1,2}:$
$$E_{a+1,2,1}^{(2)}=(E_{a+1,2,0}^{(2)})_{(E_{a+1,1,1}^{(3)})^{-1}(E_{a+1,1,0}^{(3)})^{-1}},$$
or
$$v_{a+1,2}E_{a+1,2,0}^{(2)}=(E_{a+1,1,0}^{(3)})^{-2}h_{a+1,1}^{-1}E_{a+1,2,0}^{(2)}h_{a+1,1}(E_{a+1,1,0}
^{(3)})^2$$
$$v_{a+1,2}=(E_{a+1,1,0}^{(3)})^{-2}h_{a+1,1}^{-1}(h_{a+1,1})_{(E_{a+1,2,0}^{(2)})^{-1}}(E_{a+1,1,0}^{(3)})
_{(E_{a+1,2,0}^{(2)})^{-1}}^2.$$ By abuse of notation, we remove the
last index from the $E_{.\,,\,.\,,\,.}.$
\end{proof}

We know that $\eta_{rk}$ (for $(r,k)\in I)$ is a prime element with
s.h.t. $y_{rk}$ and a central element $\mu.$ So it can be proven
(see \cite[Claim~5.5]{Te}) that
$$\eta_{r,k}=(E_{r,k-1}^{(3)})^{2}(E_{r,k-1}^{(3)})^{-2}_{(E_{r,k}^{(2)})^{-1}}$$
or
\begin{equation}\label{4.1}
\mu\eta_{r,k}^{-1}=(E_{r,k-1}^{(3)})^{-2}(E_{r,k-1}^{(3)})_{(E_{r,k}^{(2)})^{-1}}^2.\end{equation}
So we have\newpage

$$v_{r,k}\quad=(E_{r,k-1}^{(3)})^{-2}h_{r,k-1}^{-1}(h_{r,k-1})_{(E_{r,k}^{(2)})^{-1}}(E_{r,k-1}^{(3)})^2_{(E_{r,k}
^{(2)})^{-1}}$$$$
  \overset {\text{\cite[IV.6.1]{BGT5}}}=
  h_{r,k-1}^{-1}(E_{r,k-1}^{(3)})^{-2}(E_{r,k-1}^{(3)})_{(E_{r,k}^{(2)})^{-1}}(h_{r,k-1})_{(E_{r,k}^{(2)})^{-1}}$$
\begin{equation}\label{4.2} \quad=h_{r,k-1}^{-1}\mu\eta_{r,k}^{-1}(h_{r,k-1})_{(E_{r,k}^{(2)})^{-1}}.\end{equation}
We compute now $(h_{r,k-1})_{(E_{r,k}^{(2)})^{-1}}.$ We know that
$(c_{r,k-1})_{(E_{r,k}^{(2)})^{-1}}=c_{r,k-1}b_{r,k}$
(\cite[IV.6.3]{BGT5}) and
$(\zeta_{r,k-1})_{(E_{r,k}^{(2)})^{-1}}=\zeta_{r,k-1}\eta_{r,k}.$ So
$$
(h_{r,k-1})_{(E_{r,k}^{(2)})^{-1}}\overset{\text{Proposition
\ref{prop4.3}}}=(\mu\nu)^{\lambda(k-1)}(c_{r,k-1}^{k-1}\zeta_{r,k-1}^{-k+2})_{(E_{r,k}^{(2)})^{-1}}\\
$$$$\qquad=
(\mu\nu)^{\lambda(k-1)}(c_{r,k-1})_{(E_{r,k}^{(2)})^{-1}}^{k-1}(\zeta_{r,k-1})_{(E_{r,k}^{(2)})^{-1}}^{-k+2}\\
$$$$\qquad=
(\mu\nu)^{\lambda(k-1)}(c_{r,k-1}b_{r,k})^{k-1}(\zeta_{r,k-1}\eta_{r,k})^{-k+2}\\
$$\begin{equation}\label{4.3}\qquad\quad\quad\,\,\,\,\,=
(\mu\nu)^{\lambda(k-1)}\nu^{\lambda(k-1)}c_{r,k-1}^{k-1}b_{r,k}^{k-1}\mu^{\lambda(k-2)}\eta_{r,k}^{-k+2}\zeta
_{r,k-1}^{-k+2}.\end{equation}

We substitute the expressions we found in \ref{prop4.3},
\eqref{4.1}, \eqref{4.3} in \eqref{4.2}, and we get\newline (for
$(r,k)\in I$):
\begin{align*}
v_{r,k}&=(\mu\nu)^{\lambda(k-1)}\zeta_{r,k-1}^{k-2}c_{r,k-1}^{1-k}\cdot\mu\nu_{r,k}^{-1}(\mu\nu)^{\lambda(k-1)}
\nu^{\lambda(k-1)}c_{r,k-1}^{k-1}b_{r,k}^{k-1}\\
&\quad  \cdot\mu^{\lambda(k-2)}\eta_{r,k}^{-k+2}\zeta_{r,k-1}^{-k+2}\\
\begin{pmatrix} \mu^2=1\\
\nu^2=1\end{pmatrix}&=\mu^{\lambda(k-2)+1}\nu^{\lambda(k-1)}\zeta_{r,k-1}^{k-2}c_{r,k-1}^{1-k}\eta_{r,k}^{-1}
c_{r,k-1}^{k-1}b_{r,k}^{k-1}\eta_{r,k}^{-k+2}\zeta_{r,k-1}^{-k+2}\\
\begin{pmatrix} [c_{r,k-1},\eta_{r,k}^{-1}]\\
=\nu\end{pmatrix}&=\nu^{k-1+\lambda(k-1)}\mu^{\lambda(k-2)+1}\zeta_{r,k-1}^{k-2}\eta_{r,k}^{-1}b_{r,k}^{k-1}
\eta_{r,k}^{-k+2}\zeta_{r,k-1}^{-k+2}\\
\begin{pmatrix} \forall a\in\mathbb Z,\ a+\lambda(a)\\
=\lambda(a+1)\end{pmatrix}&=
\nu^{\lambda(k)}\mu^{\lambda(k-2)+1}\zeta_{r,k-1}^{k-2}\eta_{r,k}^{-1}b_{r,k}^{k-1}\eta_{r,k}^{-k+2}\zeta_
{r,k-1}^{-k+2}\\
( [b_{r,k},\eta_{r,k}] =1)&
=\nu^{\lambda(k)}\mu^{\lambda(k-2)+1}\zeta_{r,k-1}^{k-2}b_{r,k}^{k-1}\eta_{r,k}^{-k+1}\zeta_{r,k-1}^{-k+2}\\
\begin{pmatrix} [\zeta_{r,k-1},\eta_{r,k}]=\nu\\
[\zeta_{r,k-1},b_{r,k}]=\nu\\ \mu^2=\nu^2=1\\ \forall k,
(k-2)(k-1)\equiv0\!\!\pmod2\end{pmatrix}&=\nu^{\lambda(k)}\mu^{\lambda(k-2)+1}b_{r,k}^{k-1}\eta_{r,k}^{-k+1}.
\end{align*}

Now assume that $a<r<b+a,\,\, r-a+1\le k\le b.$ Denote $k'=r-a+1.$
So by using $v_{r,k+1}=(v_{r,k})_{\tilde x_{rk}^{-1}\tilde
z_{rk}\tilde z_{r+1,k}^{-1}\tilde x_{r+1,k+1}},$ we see that
\begin{equation}\label{4.4}
v_{rk}=\mu^{\lambda(k'-2)+1}\nu^{\lambda(k')}b_{r,k}^{k'-1}\eta_{r,k}^{-k'+1}.\end{equation}

\begin{remark}\label{rem4.5}
(1) As in \cite{Mo}, we can consider a 3-point $\omega_{r,b},$\
$r<a$ and see that
$$b_{11}^b\eta_{11}^{2-b}=(\mu\nu)^{\lambda(b+1)}\mu$$
(see \cite[Proposition~11,(1)]{Mo}).

(2) If $b$ is odd, then $\mu=\nu$ \cite[Proposition~11,(3)]{Mo}.
\end{remark}

Consider now the 3-point $\omega_{a,b}.$ We know that
 $v_{a,b}=1,$ but
 $d_{a,b}=a_{a,b}^{a-b}\xi_{a,b}^{b-a}(\mu\nu)^{\lambda(b-a)}.$

 \begin{prop}\label{prop4.5} \begin{enumerate}\item[(1)] If $a\ne b,$
 then
 $\mu(\mu\nu)^{\lambda(b-a+1)}=(b_{1,1}\eta_{1,1}^{-1})^{a-b}\eta_{1,1}^{-1}.$
 \item[(2)] If $a=b,$ then $\mu=\eta_{1,1}=1.$\end{enumerate}
 \end{prop}

 \begin{proof} (1) By Remark \ref{rem4.4},
\begin{align}
1\ &=\tilde
x_{a,b}^{-2}(a_{a,b}^{a-b}\xi_{a,b}^{b-a}(\mu\nu)^{\lambda(b-a)})_{y_{a,b}^{-1}x_{a,b}^{-1}}(\tilde
x_{a,b}^2)_{y_{a,b}^{-1}}\label{4.5}\\
 \begin{pmatrix}  \tilde x_{a,b}^{-2}(\tilde x_{a,b}^2)_{y_{a,b}^{-1}}\\
 =\mu\eta_{a,b}^{-1};\\ b_{a,b}=(a_{a,b})_{y_{a,b}^{-1}x_{a,b}^{-1}}\\
 \eta_{a,b}=(\xi_{a,b})_{y_{ab}^{-1}x_{a,b}^{-1}}
 \end{pmatrix}&=b_{a,b}^{a-b}\eta_{a,b}^{b-a}\nu^{a-b}\mu^{b-a}(\mu\nu)^{\lambda(b-a)}\mu\eta_{a,b}^{-1}\nonumber\\
 &=b_{a,b}^{a-b}\eta_{a,b}^{b-a-1}(\mu\nu)^{b-a}(\mu\nu)^{\lambda(b-a)}\mu\nonumber\\
 &=b_{a,b}^{a-b}\eta_{a,b}^{b-a-1}(\mu\nu)^{\lambda(b-a+1)}\mu.\nonumber
\end{align}
So $\exists\ \gamma\in\tilde B_n$ s.t. $(b_{a,b})_\gamma=b_{1,1},$\
$(\eta_{a,b})_\gamma=\eta_{1,1}.$ Applying it, we obtain what we
wanted.

(2) By Remark \ref{rem4.4}, we have
\begin{equation}\label{4.6} 1=\tilde x_{a,b}^{-2}(\tilde
x_{a,b}^2)_{y_{a,b}^{-1}}=\mu\eta_{a,b}^{-1}\end{equation} or
$$\mu=\eta_{a,b}.$$
By the same argument as in (1), $\mu=\eta_{1,1}.$ By \eqref{4.6}, we
see that $\tilde x_{a,b}^2(\tilde x_{a,b}^{-2})_{y_{a,b^{-1}}}=1,$
or $\eta_{a,b}=1;$ that is, $\mu=\eta_{1,1}=1.$
\end{proof}

\begin{prop}\label{prop4.6} If $a\ne b,$ then
$(b_{1,1}\eta_{1,1}^{-1})^{a-b}\eta_{1,1}^{-1}=(\mu\nu)^{\lambda(b-a)}.$
\end{prop}

\begin{proof} By \eqref{4.5},
$$
\eta_{a,b}=(a_{a,b}^{a-b}\xi_{a,b}^{b-a}(\mu\nu)^{\lambda(b-a)})_{y_{a,b}^{-1}x_{a,b}^{-1}}=b_{a,b}^{a-b}
\eta_{a,b}^{b-a}(\mu\nu)^{\lambda(b-a)}.$$ Applying $\gamma$ from
above, we are done.
\end{proof}

Note that if $a=b$, we get $b_{1,1}^b=\nu^{\lambda(b+1)}.$
\begin{prop}\label{prop4.7} If $b$ is even, $a$ is odd, then $\nu=1;$ otherwise
$\mu=\nu=1.$
\end{prop}
\begin{proof} We will first prove a lemma.

\begin{lem}\label{lem4.3} $\forall\ r,\ a<r<a+b,$ we have
$$\eta_{1,1}^{b-a-1}b_{1,1}^{a-b}=\mu^{\lambda(r-a-1)}\nu^{\lambda(r-a+1)}(\mu\nu)^{\lambda(b-r+1)}.$$
\end{lem}
\begin{proof} By Remark \ref{rem4.4} and \eqref{4.4}, we have from the 3-point
$\omega_{rb}$ $(k'=r-a+1)$:
\begin{align*}\mu^{\lambda(k'-2)+1}\nu^{\lambda(k')}b_{r,b}^{k'-1}\eta_{r,b}^{-k'+1}&=\tilde
 x_{r,b}^{-2}(a_{r,b}^{r-b}\xi_{r,b}^{b-r}(\mu\nu)^{\lambda(b-r)})_{y_{r,b}^{-1}
 x_{r,b}}(\tilde x_{r,b}^2)_{y_{r,b}^{-1}}\\
 &=b_{r,b}^{r-b}\eta_{r,b}^{b-r}\nu^{r-b}\mu^{b-r}(\mu\nu)^{\lambda(b-r)}\mu\eta_{r,b}^{-1}
\end{align*}
\begin{align*} &\quad\Rightarrow
\mu^{\lambda(k'-2)}\nu^{\lambda(k')}b_{r,b}^{k'-1}\eta_{r,b}^{-k'+1}=\nu^{r-b}\mu^{b-r}(\mu\nu)^{\lambda(
b-r)}\eta_{r,b}^{b-r-1}b_{r,b}^{r-b}\\
&\quad\Rightarrow
\eta_{r,b}^{b-r-1+k'-1}b_{r,b}^{r-b-k'+1}=\mu^{\lambda(k'-2)+r-b}\nu^{\lambda(k')+b-r}(\mu\nu)^{\lambda(b-r)}\\
&\overset{k'=r-a+1}\Rightarrow
\eta_{r,b}^{b-a-1}b_{r,b}^{a-b}=\mu^{\lambda(r-a-1)}\nu^{\lambda(r-a+1)}(\mu\nu)^{\lambda(b-r+1)}
\end{align*}
$\forall\ r,$\ $\exists \gamma_r\in\tilde B_n,$ s.t.
$(\eta_{r,b})_{\gamma_r}=\eta_{1,1},$\
$(b_{r,b})_{\gamma_r}=b_{1,1}.$  Apply it, and we are done.
\end{proof}

Assume $b$ is odd. So we know that $\mu=\nu$ (By Remark
\ref{rem4.5}). If $a=b,$ then $\mu=\nu=1$ (by Proposition 4.6) Else,
$a\ne b.$ So from \lemref{lem4.3}, set $r=a+1,$ and we get
$\eta(b_{1,1}\eta_{1,1}^{-1})^{a-b}\eta_{1,1}^{-1}=\mu$. From
Proposition \ref{prop4.6}, if $\mu=\nu,$\
$(b_{1,1}\eta_{1,1}^{-1})^{a-b}\eta_{1,1}^{-1}=1,$ so
$\mu=1\Rightarrow\mu=\nu=1.$

Assume now that $b$ is even. From \lemref{lem4.3}, when setting
$r=a+1,$ we get
\begin{equation}\label{4.7}(b_{1,1}\eta_{1,1}^{-1})^{a-b}\eta_{1,1}^{-1}=\nu(\mu\nu)^{\lambda(b-a)}.\end{equation}
If $a=b,$ then $\eta_{1,1}^{-1}=\nu;$ but $\eta_{1,1}=1,$ so
$\mu=\nu=1.$ Else $(a\ne b),$ we have by Proposition \ref{prop4.6},
$(b_{1,1}\eta_{1,1}^{-1})^{a-b}\eta_{1,1}^{-1}=(\mu\nu)^{\lambda(b-a)}.$
So we have $\nu=1$ when $b$ is even. Assume now that $a$ is also
even (and $a\ne b).$ By Proposition \ref{prop4.5}, we get
\begin{equation}\label{4.8}(b_{1,1}\eta_{1,1}^{-1})^{a-b}=\mu(\mu\nu)^{\lambda(b-a+1)}\end{equation}
or (substituting $\nu=1)$, we have the set of equations:
$$\begin{cases}
(b_{1,1}\eta_{1,1}^{-1})^{a-b}\eta_{1,1}^{-1}=\mu\cdot\mu^{\lambda(b-a+1)}\\
(b_{1,1}\eta_{1,1}^{-1})^{a-b}\eta_{1,1}^{-1}=\mu^{\lambda(b-a)}\end{cases}$$
Thus, $$\mu\cdot\mu^{\lambda(b-a+1)}=\mu^{\lambda(b-a)}\Rightarrow
\mu\cdot\mu^{b-a}=1\overset{\ b-a\text{\ is
even}}\Rightarrow\mu=1.$$
\end{proof}

As in \cite{Mo}, we define a $\tilde B_n$-group $G_0(n)$ as the
subgroup of $G(n)$ generated by $u_1,\dots,u_{n-1};$\ $G_0(n)$ is
$\tilde B_n$-isomorphic to $\tilde P_{n,0}$ (recall that
$n=2ab+b^2).$

\begin{defn}\label{defn4.11} $G_0(n)$ is a group with

\noindent Generators: \begin{align*}M_0=&\{A_{ij},\ 1\le j\le b,\ 1\le i\le
a+j;\quad B_{ij},\ 1\le j\le b,\ 1\le i\le a+j-1;\\
&\quad C_{ij},\ a<i<a+b,\ j=i-a\}.\end{align*}

Relations:

\begin{enumerate}\item[(1)] $\forall\ a,b\in M_0$ which are adjacent,
$[a,b]=\tau$, where $\tau$ is independent of (such) $a,b$,\
$\tau^2=1,$\ $\tau_d=\tau\ \forall\ d\in M_0.$ \item[(2)] If $a,b\in
M_0$ are not adjacent, then $[a,b]=1$.
\end{enumerate}

 for each $d\in M_0$ we
introduce the notion of supporting half-twist from $\tilde B_n$
(resp. $B_n$) as follows: for $d=A_{ij}$, it will be $x_{ij}$ (resp.
$X_{ij})$; for $d=B_{ij},$ it will be $y_{ij}$ (resp. $Y_{ij});$ for
$d=C_{ij},$\ $i=j+a,$ it will be $z_{ij}$ (resp. $Z_{ij}).$

We say that $a,b\in M_0$ are adjacent if their supporting half-twists are
adjacent.

The $\tilde B_n$-action on $G_0(n)$ in terms of $\tilde
M=\{x_{ij},y_{ij}\}\cup\left\{z_{ij}\bigm|_{j=i-a}^{a<i<a+b}\right\}$ and $M_0$
is defined in \cite[Remark~6]{Mo}.  We consider the elements of
$$ \tilde M_1=  \tilde M\cup\{ z_{ij}|(i,j)\in\Vertices(K(a,b)),\ i,j\ge 1\
 \text{and if}\ a<i<a+b,\ \text{then}\ j\ne i-a\}$$

as polarized half-twists, and define a larger subset of $G_0(n):
\hat{M}_0 ;\,$when $ M_0 \subset \hat{M}_0$ s.t.:
$$ \hat M_0=   M_0\cup\{ C_{ij}|(i,j)\in\Vertices(K(a,b)),\ i,j\ge 1\
 \text{and if}\ a<i<a+b,\ \text{then}\ j\ne i-a\}.$$

 We start with the pair $\{B_{1,1},y_{1,1}\}$. Then $\forall t\in\tilde M_1$,
 define $L_0(t)\in\tilde M_0$ as the unique element
 $L_{\{B_{1,1},y_{1,1}\}}(t)$ s.t. $\{B_{1,1},y_{1,1}\}$ and
 $\{L_{\{B_{1,1},y_{1,1}\}}(t),t\}$ are coherent. The definition of a $\tilde
 B_n$-action on $G_0(n)$ is such that $L_0(x_{ij})=A_{ij},$\
 $L_0(y_{ij})=B_{ij},$\ $L_0(z_{ij})=C_{ij}$ where $a<i<a+b,$\ $j=i-a.$ So for
 $t\in\tilde M$, we have $L(t)\in M_0.$
 \end{defn}

 Define $C_{ij}=L_0(z_{ij}).$

 \begin{defn}\label{defn4.12} Using the $\tilde B_n$-action on $G_0(n),$ we
 define canonically the semi-direct product $G_0(n)\rtimes\tilde B_n.$ Let
 $u=y_{1,1}^2x_{2,1}^{-2}\in\tilde P_{n,0}\subset\tilde B_n.$ Let $N(a,b)$ be the
 normal subgroup of $G_0(n)\rtimes \tilde B_n,$ normally generated by the
 elements:
 \begin{align*}
 &n_1=B_{1,1}^bu^{2-b}c(c\tau)^{\lambda(b+1)};\\
 &n_2=(c\tau)^b;\\
 &n_3=(B_{1,1}u^{-1})^{a-b}u^{-1}\cdot c(c\tau)^{\lambda(b-a+1)};\\
 &n_4=(B_{1,1}u^{-1})^{a-b}u^{-1}\cdot \tau(c\tau)^{\lambda(b-a)}
 \end{align*}
 (when $c=[x^2,y^2],$ \ $x,y$ are any two adjacent half-twists in $\tilde
 B_n;$ \ $\lambda(k)=\frac{k(k-1)}{2}).$

 Note that the elements in $N(a,b)$ are defined according to the
 relations found in Proposition 4.8 (\eqref{4.7}, \eqref{4.8}) and
 Remark 4.5.

  \end{defn}

 Define $$G(a,b)=(G_0(n)\rtimes\tilde B_n)/N(a,b).$$
 So as in \cite[Proposition~32]{Mo}, one can prove that
 $$\pi_1(\C^2-S_{F_{1,(a,b)}})\simeq G(a,b).$$

 Define $\psi_{a,b}:G(a,b)\to S_n,$ by $\psi_{a,b}(\alpha,\beta)=\psi(\beta)$
 where $\psi: \tilde B_n\to S_n$ is the homomorphism to the symmetric group,
 induced from the standard homomorphism $B_n\to S_n$. Let
 $H_{a,b}=\ker\psi_{a,b},$\ $(H_{a,b})_0=\ker(H_{a,b}\to Ab(G(a,b))),$ or, in
 other words, if $Ab_{a,b}=$ abelization map of $G(a,b),$ then
 $(H_{a,b})_0=\ker\psi_{a,b}\cap \ker Ab_{a,b}.$ Note that
 $G(a,b)/H_{a,b}\simeq S_n.$
 Also define $\bar{\psi}_{a,b}:\pi_1(\C\P^2 -
 \bar{S}_{F_{1,(a,b)}})\to S_n$, and let
 $\overline{H}_{a,b}=\ker\bar{\psi}_{a,b}$. In the same way as above, we
 define $(\overline{H}_{a,b})_0$ and $(\overline{H}_{a,b})'_0$.

  So we have the following

 \begin{thm}\label{4.1} {}\quad

 \begin{enumerate}\item[1)] $H_{a,b}/(H_{a,b})_0\simeq\Z.$
\item[2)] $H_{a,b}'=(H_{a,b})_0'\simeq\begin{cases} \Z_2\  & b\ \text{even},\
a\ \text{odd}\\
1\  & \text{else}\end{cases}$\newline $H_{a,b}'\subset\Center(G(a,b)).$
\item[3)] $Ab(H _{a,b})_0\simeq(\Z_{b-2a})^{n-1}.$ \end{enumerate}
\end{thm}

\begin{proof}  The statement can be deduced directly from the definition of
$G(a,b).$\ 2) follows from Proposition \ref{prop4.7}. 3) follows
from the definition of $N(a,b)$ and the following facts:\\
$n_1= B_{1,1}^bu^{2-b}c(c\tau)^{\lambda(b+1)} =
(B_{1,1}u^{-1})^bu^2c(c\tau)^{\lambda(b+1)}$ and
\begin{align*}\Z^2/\langle(b,2),(a-b,-1)\rangle&=\Z^2/\langle
(b,2),(a,-1)\rangle=\Z^2/\langle(b-2a,0),(a-b,1)\rangle\\
&=\Z^2/\langle(b-2a,0),(0,1)\rangle=\Z_{b-2a},\end{align*}
\end{proof}
As in \cite[p.~74]{Mo}, one can consider the projective case
$$\pi_1(\C\P^2-\bar{S}_{F_{1,(a,b)}})\simeq G(a,b)/(y_{1,0}^{2m_1}\cdot U),$$
where $2m_1=\deg \bar{S}_{F_{1,(a,b)}}=6ab-2a-2b-3b+3b^2,$\ $U\in
(H_{a,b})_0.$ From the definition of $\overline{H}_{a,b}$,
$(\overline{H}_{a,b})_0$ it follows that they coincide with the
images of $H_{a,b}$ and $(H_{a,b})_0$ in
$G(a,b)/(y_{1,0}^{2m_1}\cdot U) = \overline G(a,b)$. So by the same
arguments as in \cite{Mo}, we have
$$\overline H_{a,b}/(\overline{H_{a,b}})_0\simeq\Z_{m_1},\quad (\overline
H_{a,b})_0\simeq (H_{a,b})_0,$$ so
$$Ab(\overline H_{a,b})_0\simeq (\Z_{b-2a})^{n-1},$$
and
$$\overline H_{a,b}'\simeq (\overline H_{a,b})_0'\simeq
(H_{a,b})'_0\simeq\begin{cases} \Z_2\ &b\ \text{even} ,\
a\ \text{odd} \\
1\ & \text{else}\end{cases}$$ Thus, there exists  a series
$$1\vartriangleleft(H_{a,b})_0'\vartriangleleft (H_{a,b})_0\vartriangleleft
H_{a,b}\vartriangleleft G(a,b)$$ s.t.
 \begin{align*}
&G(a,b)/H_{a,b}\simeq S_n\\
&H_{a,b}/(H_{a,b})_0\simeq \Z\\
&(H_{a,b})_0/(H_{a,b})_0'\simeq(\Z_{b-2a})^{n-1},\end{align*} and
$$(H_{a,b})_0'\simeq\begin{cases}\Z_2\ &b\ \text{even},\ a\ \text{odd}\\
1\ &\text{else}
\end{cases}$$

and a series:
$$1\vartriangleleft(\overline H_{a,b})_0'\vartriangleleft (\overline H_{a,b})_0\vartriangleleft
\overline H_{a,b}\vartriangleleft \overline G(a,b)$$ s.t.

\begin{align*}
&\overline G(a,b)/ \overline H_{a,b}\simeq G(a,b)/H_{a,b}\\
&\overline H_{a,b}/ ( \overline H_{a,b})_0\simeq \Z_{m_1}\\
&(\overline H_{a,b})_0/(\overline H_{a,b})_0' \simeq ( H_{a,b})_0/(
H_{a,b})_0',\end{align*} and
$$(\overline H_{a,b})_0' \simeq ( H_{a,b})_0'$$

\section{Appendix}
This Appendix describes the braid monodromy factorizations induced
from the regeneration of each point and the induced relations from
it.

For computing explicitly the braid monodromy factorizations
$\mathcal{H}(r,k)$ induced from the
 6/3/2-points - $\omega_{r,k}\,$, we use the results of \cite{Mo}.

 For $(r,k) = (0,0),(a,0)$, the vertex $\omega_{r,k}$ is a 2--point on
 the edge $L_j$ (a point which is on the intersection of two
 planes). Therefore, the braid monodromy factorization of the
 regenerated neighborhood of the vertex $\omega_{r,k}$ is $$\mathcal{H}(r,k) =
 Z_{j,j`}.$$

 For $(r,k)$ s.t $\omega_{r,k}\,$ are on the boundary of $P$ and
 $(r,k) \neq (0,b),(a+b,b),(0,0),(a,0)$ -- $\omega_{r,k}$ is a 3--point
 (a point that lies on the intersection of three
 planes), such that locally it looks like one of the following configurations:
\begin{center}
 \epsfig{file=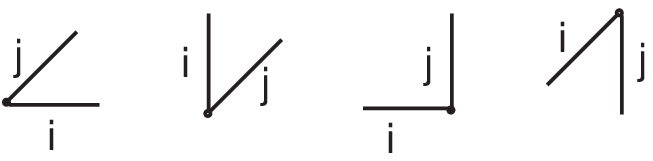}
\end{center}
Consider the first and the third cases (where the line $L_j$ is
regenerated first). Then the braid monodromy factorization of the
 regenerated neighborhood of the vertex $\omega_{r,k}$ is $$\mathcal{H}(r,k)
 = Z^{(3)}_{i\,i`,j}\tilde{Z}_{j\,j`(i)}$$
 when $Z^{(3)}_{i\,i`,j} =
 Z^3_{i`,j}Z^3_{i,j}(Z^3_{i,j})_{Z_{i,i'}}$.

 Consider the second and the fourth cases (where the line $L_i$ is
regenerated first). Then the braid monodromy factorization of the
 regenerated neighborhood of the vertex $\omega_{r,k}$ is
 $$\mathcal{H}(r,k) = Z^{(3)}_{j\,j`,i}\tilde{Z}_{i\,i`(j)}.$$

 In both cases, $\tilde{Z}_{j\,j`(i)}$ is represented by the following path:
 \begin{center}
 \epsfig{file=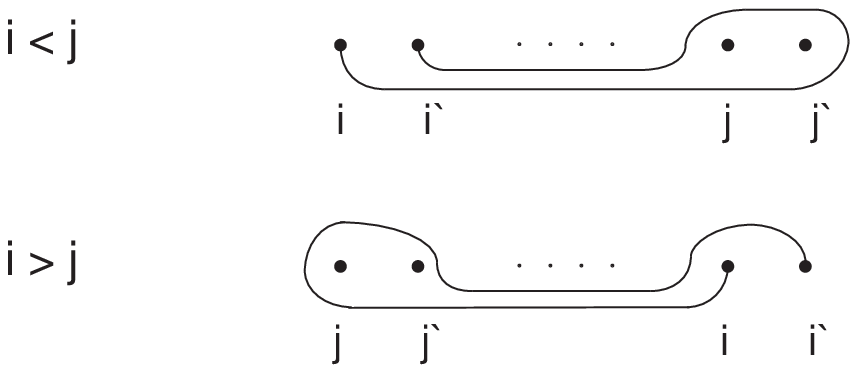}
\end{center}

For $(r,k)$ such that $\omega_{r,k}\,$ are not on the boundary of
$P$, $\omega_{r,k}$ is a 6--point. Assume that locally it looks like
the following configuration (when the lines are numerated locally):
 \begin{center}
 \epsfig{file=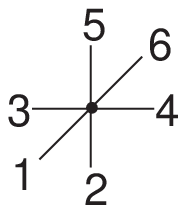}
\end{center}

Then the braid monodromy factorization of the
 regenerated neighborhood of the vertex $\omega_{r,k}$ is:

\begin{center}
$\mathcal{H}(r,k)=Z^{(3)}_{1',2\,2'}\tilde{Z}_{6\,6'}Z^{(2)}_{3\,3',6'}(Z^{(2)}_{2\,2',6'})^{\bullet}
\bar{Z}^{(3)}_{4\,4',6}(Z^{(2)}_{3\,3',6})^{\bullet}(Z^{(2)}_{2\,2',6})^{\bullet}
(\hat{F}(\hat{F})_{\rho^{-1}}))^{\bullet}Z^{(3)}_{5\:5',6}$\\
 $\Bigg(\prod\limits_{i=6',6,5'
\atop5,4',4}(Z^2_{1',i})\Bigg)^{\bullet}
\,\bar{Z}^{(3)}_{1',3\,3'}\prod\limits_{i=6',6,5'
\atop5,4',4}(Z^2_{1\,i})\tilde{Z}_{1,1'},$\end{center} where
$Z^{(2)}_{i\,i`,j} = Z^2_{i`,j}Z^2_{i,j}\,$,$( )^{\bullet}$ is the
conjugation by the braid induced from the motion:
\begin{center}
\epsfig{file=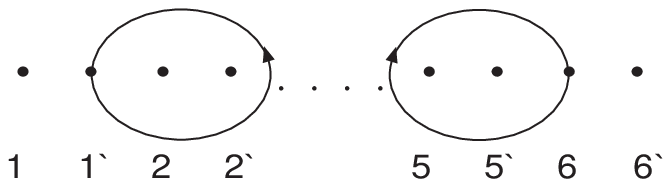}\\
\end{center}
and $\tilde{Z}_{1\,1'}, \tilde{Z}_{6\:6'}$ are
\begin{center}
\epsfig{file=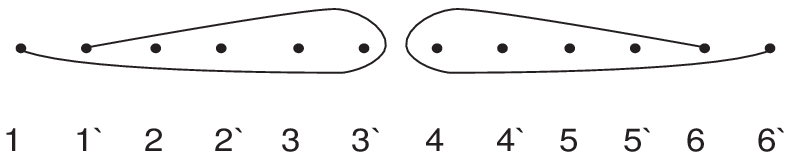}\\
\end{center}
$\rho = Z_{2\,2'}Z_{5\:5'}$\\
$\hat{F}=Z^{(3)}_{2',3\,3'}Z^{(3)}_{4\,4',5}\check{Z}_{3'\,4}\check{Z}_{3\,4'}
\overset{(3-3')}{Z^2_{2',5}}\bar{Z}^2_{2',5'}$\\
where $\check{Z}_{3\,4'},\,\check{Z}_{3'\,4}$ are:
\begin{center}
\epsfig{file=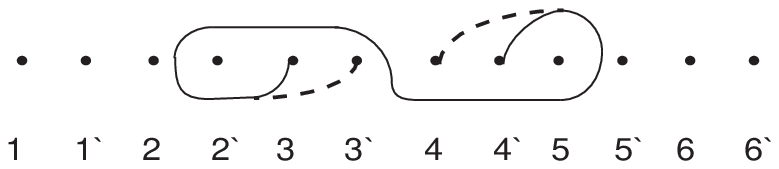}\\
\end{center}

By the Van-Kampen Theorem \cite{vK}, we can see that we get a triple
relation ($\langle A,B \rangle = e$) for each pair of generators
whose corresponding lines (from which they are created) induce a
common triangle in the complex $K(a,b)$; and we get a double
(commutation) relation ($[A,B] = e$) for each pair of generators
whose corresponding lines does not induce a common triangle in the
complex. This is the basis for the embedding of $\tilde{B_n}$ in
$G$. For more details, see \cite{Mo}.

\bigskip
\textsc{Michael Friedman, Department of Mathematics, Bar-Ilan
University, 52900 Ramat Gan, Israel\\}
\textsl{email}: fridmam@macs.biu.ac.il \\\\
\textsc{Mina Teicher, Department of Mathematics, Bar-Ilan
University, 52900 Ramat Gan, Israel\\} \textsl{email}:
teicher@macs.biu.ac.il

\end{document}